\documentclass{article}
\usepackage{amssymb}
\newtheorem{corollary}{Corollary}
\newtheorem{theorem}{Theorem}
\newtheorem{proposition}{Proposition}
\usepackage{color}
\usepackage{graphicx}
\usepackage{float}

\begin{document}

\begin{center} {\large Some notes on the pseudorandomness of Legendre symbol and Liouville function}
\end{center}

\begin{center} Johannes Grünberger and Arne Winterhof\\
Johann Radon Institute for Computational and Applied Mathematics\\
Austrian Academy of Sciences\\
Linz, Austria\\
Email: arne.winterhof@ricam.oeaw.ac.at
\end{center}

\begin{abstract}
We improve bounds on the degree and sparsity of Boolean functions representing the Legendre symbol as well as on the $N$th linear complexity of the Legendre sequence. We also prove similar results for both the Liouville function for integers and its analog for polynomials over $\mathbb{F}_2$, or more general for any (binary) arithmetic function which satisfies 
$f(2n)=-f(n)$ for $n=1,2,\ldots$
\end{abstract}

Keywords: arithmetic functions, Legendre symbol, Liouville function, Boolean functions, sequences, algebraic degree, algebraic thickness, linear complexity, lattice test\\

AMS-class: 11T71, 06E30, 94A55\\

\section{Introduction}

Several (binary) {\em arithmetic functions}, that is, mappings from $\mathbb{N}$ to $\{-1,+1\}$ possess some properties of pseudorandomness, that is, they are not distinguishable from a truly random function with respect to certain measures of pseudorandomness such as balance and correlation. 

In this paper we consider the Legendre symbol and
the Liouville function for integers and polynomials, respectively. We will recall the definitions and summarize some known features of pseudorandomness of these functions in Section~\ref{survey}.

Note that for $r=1,2,\ldots$ the first $2^r-1$ values of a (binary) arithmetic function $f:\mathbb{N}\rightarrow \{-1,+1\}$ can be identified with a Boolean function $B:\mathbb{F}_2^r\rightarrow \mathbb{F}_2$ in $r$ variables,
 via the equation
$$f\left(\sum_{j=1}^r n_j 2^{j-1}\right)=(-1)^{B(n_1,n_2,\ldots,n_r)},\quad (n_1,n_2,\ldots,n_r)\in \mathbb{F}_2^r \setminus \{(0,0,\ldots,0)\},$$
which is unique up to the value $B(0,0,\ldots,0)$.
Here we identify the finite field $\mathbb{F}_2$ of two elements with the set of integers $\{0,1\}$.
Moreover, every Boolean function in $r$ variables can be identified with a unique polynomial over $\mathbb{F}_2$ in $r$ variables with all local degrees at most $1$,
$$B(X_1,X_2,\ldots,X_r)=\sum_{I\subseteq \{1,2,\ldots,r\}} a_I X^I,$$
where
$$X^I=\prod_{j\in I} X_j$$
and $a_I\in \mathbb{F}_2$.
The {\em algebraic degree} of a nonzero Boolean function $B$ is 
$$\deg(B)=\max\{|I| : I \subseteq \{1,2,\ldots,r\} \mbox{ with }a_I=1\}$$
with the convention $\deg(0)=0$,
and the {\em sparsity} (or {\em algebraic thickness}) of $B$ is the number of nonzero coefficients of $B$,
$${\rm spr}(B)=|\{I \subseteq \{1,2,\ldots,r\}: a_I=1\}|.$$
The expected value of the sparsity is
$$\frac{1}{2^{2^r}}\sum_B {\rm spr}(B)=\frac{1}{2^{2^r}}\sum_{s=0}^{2^r}s{2^r\choose s}=2^{r-1},\quad r=1,2,\ldots$$
where the first sum is over all Boolean functions in $r$ variables.
For the expected value of the algebraic degree we have
\begin{eqnarray*}\frac{1}{2^{2^r}}\sum_B \deg(B)&=&\frac{1}{2^{2^r}} \sum_{d=1}^r d \left(2^{{r\choose d}}-1\right)2^{\sum_{j=0}^{d-1}{r\choose j}}\\
&\ge& \frac{1}{2^{2^r}}\left((r-1)(2^r-1)2^{2^r-r-1}+r2^{2^r-1}\right)\\
&=&r-\frac{1}{2}-\frac{r-1}{2^{r+1}} \ge r-\frac{5}{8}.
\end{eqnarray*}
Thus a pseudorandom Boolean function must have both large sparsity and large degree. Both are also required for Boolean functions used in cryptography, see the monographs \cite{C20,CS17,M16,S03}. 

In Section~\ref{bool} we will observe that the above mentioned arithmetic functions satisfy these desirable features of pseudorandomness. In particular, we improve some results for the Legendre symbol of \cite[Chapter~10]{S03} in the case that the least quadratic non-residue ${\rm N}(p)$ modulo $p$ is small, say, ${\rm N}(p)\le 17$.
The least quadratic non-residue is obviously a prime and from Dirichlet's theorem and the law of quadratic reciprocity we get
$$\lim_{x\rightarrow \infty}\frac{\left|\{p\le x: {\rm N}(p)=p_k\}\right|}{\pi(x)}=\frac{1}{2^k},\quad k=1,2,\ldots,$$
where 
$p_k$ denotes the $k$th prime and $\pi(x)$ is the number of primes at most $x$, see \cite{MT21}.
For example, we have 
$$\lim_{x\rightarrow \infty}\frac{\left|\{p\le x: {\rm N}(p)\le 17\}\right|}{\pi(x)}=\frac{127}{128}>0.99$$
and the primes $p$ with ${\rm N}(p)>17$ are quite rare compared to the primes with ${\rm N}(p)\le 17$.
Recall that ${\rm N}(p)=2$ if and only if $p\equiv \pm 1\bmod 8$.

We can also identify (binary) arithmetic functions $f:\mathbb{N}\rightarrow \{-1,+1\}$ with (binary) sequences $(s_n)_{n=1}^\infty$, $s_n\in \mathbb{F}_2$, $n=1,2,\ldots$, via
\begin{equation}\label{funcseq}
f(n)=(-1)^{s_n},\quad n=1,2,\ldots
\end{equation}

For $N=1,2,\ldots$ the {\em $N$th linear complexity} $L({\cal S},N)$ of a binary sequence ${\cal S}=(s_n)_{n=1}^\infty$
is the smallest positive integer $L$ such that there are constants $c_0,c_1,\ldots,c_{L-1}\in \mathbb{F}_2$
with
$$s_{n+L}=c_{L-1}s_{n+L-1}+\ldots+c_0s_n,\quad n=1,2,\ldots,N-L,$$
with the convention $L({\cal S},N)=0$ if $s_1=s_2=\ldots=s_n=0$ and $L({\cal S},N)=N$ if $s_1=s_2=\ldots=s_{N-1}=0$ and $s_N=1$.
The {\em linear complexity} $L({\cal S})$ is
$$L({\cal S})=\sup_{N=1,2,\ldots} L({\cal S},N).$$
If ${\cal S}$ is $T$-periodic, we have $L({\cal S})\le T$, and $L({\cal S})<\infty$ if and only if ${\cal S}$ is ultimately periodic.

The expected value of the $N$th linear complexity is (with a slight abuse of notation)
$$\frac{1}{2^N}\sum_{{\cal S}\in \mathbb{F}_2^N}L({\cal S},N)=\frac{N}{2}+O(1),$$
see \cite{G76},
where 
$f(N)=O(g(N))$ is equivalent to $|f(N)|\le c g(N)$ for some $c>0$. 
(Note that here $N$ is fixed.)
Niederreiter \cite{N88} showed that a random sequence~${\cal S}$ satisfies
$$L({\cal S},N)=\frac{N}{2}+O(\log N),$$
where deviations from $\frac{N}{2}$ of order of magnitude $\log N$ must appear for infinitely many $N$. Here $\log N$ is the natural logarithm of $N$.
(Here ${\cal S}$ is fixed.)

In Section~\ref{seq} we improve a result of \cite{CGGT22} which already improved the bound on the $N$th  linear complexity of the Legendre sequence of \cite[Theorem~9.2]{S03} by a factor $\log N$. 
For $p\equiv \pm 3\bmod 8$ we show that the $N$th linear complexity of the Legendre sequence with modulus $p$ is between $(N-1)/2$ and $N/2+1$ for $1\le N\le 2p-1$
which substantially improves the lower bound of order of magnitude
$\frac{\min\{N,p\}}{p^{1/2}}$ of \cite{CGGT22}
to the best possible order of magnitude $N$.
We also prove analog results for the Liouville functions.
With respect to the result of \cite{N88}, the Legendre symbol for $p\equiv \pm 3\bmod 8$ and both Liouville functions do not quite behave like random functions.

For $p\equiv \pm 1\bmod 8$ we improve the bound of \cite{CGGT22} on the $N$th linear complexity of the Legendre sequence of modulus $p$ by a factor $\log p$
and we provide some numerical data which supports the conjecture
that also in this case the $N$th linear complexity follows closely $N/2$
but in addition its maximal deviation from $N/2$ is of order of magnitude $\log p$.
Hence, for $p\equiv \pm 1\bmod 8$ we expect that the Legendre sequence behaves like a random (periodic) sequence (of period~$p$).

\section{Some arithmetic functions and their pseudorandomness}\label{survey}

\subsection{Legendre symbol}

For a prime $p>2$ the {\em Legendre symbol} $\left(\frac{n}{p}\right)$ is defined by 
$$\left(\frac{n}{p}\right)=\left\{\begin{array}{cc} 1,& n \mbox{ is a quadratic residue modulo $p$},\\
-1,& n \mbox{ is a quadratic non-residue modulo $p$},\\
0, & n\equiv 0\bmod p.\end{array}\right.$$ 
We can identify the Legendre symbol with a binary arithmetic function $\ell$ defined by
$$\ell(n)=\left\{\begin{array}{cc} 1, & n \mbox{ is a quadratic residue modulo $p$ or $n\equiv 0\bmod p$},\\
-1, & n\mbox{ is a quadratic non-residue modulo $p$}.\end{array}\right.$$

The Legendre symbol possesses several features of pseudorandomness.

The Legendre symbol is (locally) balanced by the {\em Burgess bound}, see for example \cite[(12.58)]{IK04}
$$\sum_{n=1}^N \left(\frac{n}{p}\right)=O\left(N^{1-1/r}p^{(r+1)/(4r^2)}\log(p)^{1/r}\right),\quad N=1,\ldots,p-1,\quad r=1,2,\ldots,$$
which implies
$$\sum_{n=1}^N \left(\frac{n}{p}\right)=o(N) \quad\mbox{for }N\ge p^{1/4+o(1)},$$
where 
$$f(N)=o(g(N))\quad\mbox{if}\quad \lim_{N\rightarrow\infty}\frac{f(N)}{g(N)}=0.$$
Note that the least quadratic non-residue ${\rm N}(p)$ modulo $p$ is
$${\rm N}(p)=O\left(p^{1/(4e^{1/2})+o(1)}\right),$$
see for example \cite[p.\ 156]{FI10}.
Assuming the generalised Riemann hypothesis we have
$${\rm N}(p)=O((\log p)^2),$$
see for example \cite[Theorem 8.5.3]{BS96} or Ankeny's original paper \cite{A52}. Anyway, as mentioned in the introduction, with very high probability ${\rm N}(p)$ is very small.

The Legendre symbol is uncorrelated, that is, for a fixed positive integer $k$ we have
\begin{equation}\label{corrleg}\sum_{n=1}^N\prod_{j=1}^k \left(\frac{n+d_j}{p}\right)=O(kp^{1/2}\log p)
\end{equation}
for any integers $0\le d_1<d_2<\ldots<d_k$ and $1\le N\le p-d_k-1$, see 
\cite{MS97}.

Denote by $\log_2$ the binary logarithm.
Let $B$ be the Boolean function in $r=\lfloor \log_2 p\rfloor$ variables defined by
$$\left(\frac{\sum_{j=1}^rn_j2^{j-1}}{p}\right)=(-1)^{B(n_1,n_2,\ldots,n_r)},\quad (n_1,n_2,\ldots,n_r) \in \{0,1\}^r\setminus\{(0,0,\ldots,0)\},$$
and $B(0,0,\ldots,0)$ either $0$ or $1$.
Then 
\begin{equation}\label{sprigor}{\rm spr}(B)\ge 2^{-3/2}p^{1/4}(\log  p)^{-1/2}-1
\end{equation}
and
\begin{equation}\label{degigor}\deg(B)\ge 0.041 \log_2 p+o(\log p),
\end{equation}
see \cite[Theorem 10.1 and $(10.3)$]{S03}.
For $p\equiv \pm 3\bmod 8$, that is, ${\rm N}(p)=2$, \cite[(10.4)]{S03} provides the better bounds
$${\rm spr}(B)\ge 2^{r-2}>\frac{p}{8}\quad\mbox{and}\quad \deg(B)\ge r-1>\log_2 p-2,$$
which we improve in Corollary~\ref{corspr} below. We will also improve $(\ref{sprigor})$ for all primes with ${\rm N}(p)\le 7$ and $(\ref{degigor})$ for all primes with ${\rm N}(p)\le 31$, see Theorem~\ref{thmspr} below.

Let ${\cal L}_p=(\ell_n)_{n=1}^\infty$ be the sequence  identified with the Legendre symbol $\left(\frac{n}{p}\right)$ via $(\ref{funcseq})$ for $n\not \equiv 0\bmod p$ and $\ell_{kp}=0$ for $k=1,2,\ldots$
From \cite[Corollary~4]{CGGT22} and the bound on the correlation measure of order $k$ of \cite{MS97}, that is essentially $(\ref{corrleg})$,
we get
$$L({\cal L}_p,N)\gg \frac{\min\{N,p\}}{p^{1/2}}.$$
For more details see the Appendix 2 of this paper. Here $g(N)\gg f(N)$ is equivalent to $f(N)=O(g(N))$.\\
Note that using the bound of \cite[Theorem~3.1]{NW02} we can also get non-trival bounds on the correlation measure of order $k$ for $1\le k<Np^{-1/2}$ whereas \cite{MS97} is only non-trivial for $1\le k<Np^{-1/2}(\log p)^{-1}$. 
More precisely, we get from \cite[Theorem~3.1]{NW02},
$$\sum_{n=1}^N\prod_{j=1}^k\left(\frac{n+d_j}{p}\right)=O\left(k^{1/2}N^{1/2}p^{1/4}\right).$$
Combining this bound with \cite[Corollary 4]{CGGT22} we get
\begin{equation}\label{leglow}
L({\cal L}_p,N)\gg \frac{\min\{N,p\}\log p}{p^{1/2}},\quad N\ge p^{1/2}.
\end{equation}
For $p\equiv \pm 3\bmod 8$ we will prove the improvement
$$L({\cal L}_p,N)=\frac{N}{2}+O(1),$$
see Corollary~\ref{cor1} below.

Note that we know the exact value of the linear complexity in all cases, see~\cite{DHS98,T64}:
\begin{equation}\label{leglincompl}
L({\cal L}_p)=\left\{\begin{array}{cl}(p-1)/2, & p\equiv 1 \bmod 8,\\
p, & p\equiv 3\bmod 8,\\
p-1,&p\equiv -3\bmod 8,\\
(p+1)/2, & p\equiv -1\bmod 8.\end{array}\right.
\end{equation}

For further features of pseudorandomness of the Legendre symbol we refer to the recent survey \cite{W23}.

 Note that in most references the Legendre sequence is a shift ${\cal L}_p'$ by one position of the sequence ${\cal L}_p$ studied here. In the periodic case there is no difference and $(\ref{leglincompl})$ holds in both cases.
In the aperiodic case, it is easy to see that we have
$$L({\cal L}_p',N)\le L({\cal L}_p,N)+1\le L({\cal L}_p',N+1)$$
and all results on  $L({\cal L}_p,N)$ mentioned in this paper differ by at most $1$ from the analogical results for $L({\cal L}_p',N)$.

\subsection{Liouville function for integers}
Let 
$$n=\prod_{j=1}^s p_j^{a_j}$$ 
be the unique prime factorization of an integer $n>1$ with primes $p_1<p_2<\ldots<p_s$ and positive integers $a_1,a_2,\ldots,a_s$.
Then the {\em Liouville function} $\lambda$ of $n$ is 
$$\lambda(n)=(-1)^{\sum_{j=1}^s a_j},\quad n=2,3,\ldots$$
and $\lambda(1)=1$. 
The Liouville function possesses some properties of pseudorandomness. It is asymptotically balanced
$$\sum_{n=1}^N \lambda(n)=o(N)$$
and the Riemann hypothesis is equivalent to 
$$\sum_{n=1}^N\lambda(n)=O(N^{1/2+\varepsilon})\quad \mbox{for any }\varepsilon>0,$$
see \cite{H13}, where the implied constant depends only on $\varepsilon$.

The {\em Chowla conjecture} asserts that
$$\sum_{n=1}^{N}\lambda(n+h_1)\lambda(n+h_2)\cdots \lambda(n+h_k)=o(N)$$
for any fixed $k=1,2,\ldots$ and integers $0\le h_1<h_2<\ldots<h_k$, see \cite{C65}.

Note that 
$$\lambda(2n)=-\lambda(n),\quad n=1,2,\ldots$$

\subsection{Liouville function for polynomials}

Let $F(X)$ be a non-constant polynomial over $\mathbb{F}_2$ and
$$F(X)=\prod_{j=1}^s I_j(X)^{a_j}$$
be its unique factorisation into distinct $\mathbb{F}_2$-irreducible monic polynomials $I_1$, $I_2,\ldots,I_s$ with positive integers $a_1,a_2,\ldots,a_s$.
Then the {\em polynomial Liouville function} $\lambda$ of $F$ is
$$\lambda(F)=(-1)^{\sum_{j=1}^s a_j}$$
and $\lambda(1)=1$.
We can identify $\lambda$ with an arithmetic function $\ell$ by
$$\ell\left(\sum_{j=1}^rn_j2^{j-1}\right)=\lambda\left(\sum_{j=1}^rn_jX^{j-1}\right),\quad n_1,n_2,\ldots,n_r\in \mathbb{F}_2.$$

By Carlitz~\cite{C32} we have the following property of balance,
$$\sum_{\deg(F)=d}\lambda(F)= 2^{\lfloor (d+1)/2\rfloor}.$$
For large finite fields and polynomials of fixed degree the analog of the Chowla conjecture for the polynomial Liouville function was settled in \cite{CR14} for finite fields of odd characteristic, see also \cite{MW16},
and by \cite{C15} for finite fields of even characteristic. (Actually, \cite{C15,CR14} deal with the Möbius function but the proofs and results are exactly the same for the Liouville function.)
For finite fields of odd characteristic and fixed size of extension degree at least $3$ and polynomials of sufficiently large degree see the breakthrough paper \cite{SS22}. However, the Chowla conjecture for polynomials over $\mathbb{F}_2$ seems to be still out of reach.

We also have
$$\ell(2n)=\lambda(F_{2n})=\lambda(XF_n)=-\lambda(F_n)=-\ell(n),\quad n=1,2,\ldots$$
for the Liouville function of polynomials over $\mathbb{F}_2$, where
$$F_k(X)=\sum_{j=1}^r k_jX^{j-1}\quad\mbox{if}\quad k=\sum_{j=1}^rk_j2^{j-1}, \quad k_1,k_2,\ldots,k_r\in \{0,1\}.$$

\section{Bounds on degree and sparsity}\label{bool}
Now we prove (optimal) lower bounds on degree and sparsity for arithmetic functions with $f(2n)=-f(n)$ for $n=1,2,\ldots,2^{r-1}-1$.

\begin{theorem}\label{bool1}
Let $f$ be any (binary) arithmetic function with 
$$f(2n)=-f(n),\quad n=1,2,\ldots, 2^{r-1}-1,$$
and $B$ be the Boolean function defined by
$$f\left(\sum_{j=1}^r n_j2^{j-1}\right)=(-1)^{B(n_1,n_2,\ldots,n_r)},\quad (n_1,n_2,\ldots,n_r)\in \{0,1\}^r\setminus \{(0,0,\ldots,0)\}$$
and $B(0,0,\ldots,0)=c$ with $c\in \mathbb{F}_2$.
Then we have 
$$\deg(B)\ge r-1$$
and
$${\rm spr}(B)\ge  \left\lfloor \frac{2^r}{3}\right\rfloor.$$
\end{theorem}
Proof.
The condition $f(2n)=-f(n)$ for $n=1,2,\ldots,2^{r-1}$ is equivalent to 
$$B(n_1,n_2,\ldots,n_{r-1},0)=1-B(0,n_1,\ldots,n_{r-1})$$
for $(n_1,n_2,\ldots,n_{r-1})\in \{0,1\}^{r-1}\setminus\{(0,0,\ldots,0)\}$.
The Boolean function in $r-1$ variables defined by
$$F(X_1,X_2,\ldots,X_{r-1})=B(X_1,X_2,\ldots,X_{r-1},0)+B(0,X_1,\ldots,X_{r-1})$$
satisfies $F(0,0,\ldots,0)=0$ (which does not depend on the actual value of $B(0,0,\ldots,0)$) and 
$$F(n_1,n_2,\ldots,n_{r-1})=1,\quad (n_1,\ldots,n_{r-1})\not=(0,0,\ldots,0).$$
Hence, $F$ is uniquely represented by the polynomial of local degrees at most $1$,
$$F(X_1,X_2,\ldots,X_{r-1})=1+\prod_{i=1}^{r-1}(X_i+1)$$
of degree $r-1$ and sparsity $2^{r-1}-1$.
Thus
$$\deg(B)\ge \deg(F)=r-1.$$
Write 
$$B(X_1,X_2,\ldots,X_r)=\sum_{I\subseteq \{1,2,\ldots,r\}}a_IX^I.$$ Then we have
$$F(X_1,X_2,\ldots,X_{r-1})=\sum_{\emptyset \not= I\subseteq \{1,2,\ldots,r-1\}}(a_I+a_{I+1})X^I,$$
where $I+1=\{j+1: j\in I\}$. In particular, we have
$$a_I+a_{I+1}=1 \quad \mbox{for all }I \mbox{ with }\emptyset\not= I \subseteq \{1,2,\ldots,r-1\}.$$
For fixed $I\subseteq \{1,2,\ldots,r-1\}$ with $1\in I$
put $m_I=r-\max\{j\in I\}$. Then we have 
$$a_I=a_{I+2}=\ldots=a_{I+2\lfloor m_I/2\rfloor}\not= a_{I+1}=a_{I+3}+a_{I+2\lfloor(m_I-1)/2\rfloor +1}.$$
The number of nonzero coefficients is minimal if $a_\emptyset = 0$ and $a_I = 0$ for all $I$ with $1 \in I$, that is, $a_I = 1$ if and only if the minimum of $I \not= \emptyset$ is even. Since there are exactly $2^{r - 2k}$ sets $I \subseteq \{1, \ldots, r\}$ with $\min\{i\in I\} = 2k$, $k = 1, \ldots, \lfloor r/2 \rfloor$, we have at least
$$ 
\sum_{k=1}^{\lfloor r/2 \rfloor} {2^{r-2k}}
= \frac{1}{3} \left(2^{r} - 2^{r - 2 \lfloor r/2 \rfloor} \right)
= \left\lfloor \frac{2^r}{3} \right\rfloor
$$
nonzero coefficients. \hfill $\Box$\\


Both bounds of Theorem~\ref{bool1} are optimal. Since the value of $B(0,0,\ldots,0)$ is not fixed by $f$, Theorem~1 applies to two different Boolean functions
$$B_1(X_1,X_2,\ldots,X_r)=\sum_{I\subseteq \{1,2,\ldots,r\}}a_IX^I$$
and 
$$B_2(X_1,X_2,\ldots,X_r)=\sum_{I\subseteq\{1,2,\ldots,r\}}(a_I+1)X^I$$
satisfying 
$$B_1(0,0,\ldots,0)\not=B_2(0,\ldots,0)$$
and
$$B_1(n_1,n_2,\ldots,n_r)=B_2(n_1,n_2,\ldots,n_r),\quad (n_1,n_2,\ldots,n_r)\not=(0,0,\ldots,0).$$
If $\deg(B_1)=r$, then $a_{\{1,2,\ldots,r\}}=1$ and $\deg(B_2)=r-1$.\\
Moreover, the Boolean function
$$B(X_1,X_2,\ldots,X_r)=\sum_{\emptyset \not= I\subseteq\{1,2,\ldots,r\}\atop \min\{i\in  I\} \equiv 0\bmod 2}X^I$$
of sparsity $\left\lfloor \frac{2^r}{3}\right\rfloor$ corresponds to an arithmetic function which satisfies the conditions of Theorem~\ref{bool1}.

Theorem~\ref{bool1} covers both Liouville functions for integers and polynomials, respectively. In the case that $2$ is a quadratic non-residue modulo $p$ we get the following result for the Legendre symbol with modulus $p$.  
\begin{corollary}\label{corspr}
For a prime $p\equiv \pm 3 \bmod 8$ put $r=\lfloor \log_2 p\rfloor$.
Let $B$ be defined by
$$\left(\frac{\sum_{j=1}^rn_j2^{j-1}}{p}\right)=(-1)^{B(n_1,n_2,\ldots,n_r)},\quad (n_1,n_2,\ldots,n_r)\in \{0,1\}^r\setminus\{(0,0,\ldots,0)\}$$
and $B(0,0,\ldots,0)=c$ with $c\in \mathbb{F}_2$.
Then we have
$$\deg(B)\ge r-1$$
and 
$${\rm spr}(B)\ge \left\lfloor\frac{2^r}{3}\right\rfloor > \left\lfloor \frac{p}{6}\right\rfloor.$$
\end{corollary}

Now we consider the case that $2$ is a quadratic residue modulo $p$.
First note that there is no analog of Theorem~\ref{bool1} or Corollary~\ref{corspr} for the condition 
$$f(2n)=f(n),\quad n=1,2,\ldots,2^{r-1}-1,$$
since the constant Boolean functions as well as the Boolean function 
$$B(X_1,X_2,\ldots,X_r)=X_1+X_2+\ldots+X_r$$
of degree $1$ and sparsity $r$
correspond to a function $f$ with this property.
However, for the Legendre symbol we can use instead the property
$$\left(\frac{{\rm N}(p)n}{p}\right)=-\left(\frac{n}{p}\right),\quad n\not\equiv 0\bmod p.$$

\begin{theorem}\label{thmspr}
For a prime $p\equiv \pm 1\bmod 8$ put 
$$r= \left\lfloor \log_2 p\right\rfloor\quad\mbox{and}\quad
s=\left\lceil \log_2 {\rm N}(p)\right\rceil,$$
where ${\rm N}(p)>2$ is the least quadratic non-residue modulo $p$.
Let $B$ be a Boolean function in $r$ variables satisfying
$$\left(\frac{\sum_{j=1}^r n_j2^{j-1}}{p}\right)=(-1)^{B(n_1,n_2,\ldots,n_r)},\quad (n_1,n_2,\ldots,n_r)\in \{0,1\}^r\setminus\{(0,0,\ldots,0)\}.$$
Then we have 
$$\deg(B)\ge \left\lfloor \frac{r}{s}\right\rfloor$$
and
$${\rm spr}(B)\ge 2^{\lfloor r/s\rfloor -1}.$$
\end{theorem}
Proof.
By the definition of $s$ we have 
$2^{s-1}<{\rm N}(p)<2^s$. We write
$${\rm N}(p)=1+\sum_{j=1}^{s-2}b_j2^j+2^{s-1} \quad\mbox{with}\quad b_1,\ldots,b_{s-2}\in \{0,1\}.$$
From 
$$\left(\frac{{\rm N}(p)n}{p}\right)=-\left(\frac{n}{p}\right),\quad n=1,2,\ldots,p-1,$$
we get
\begin{eqnarray*}
&&B(n_1,\underbrace{0,\ldots,0}_{s-1},n_2,\underbrace{0,\ldots,0}_{s-1},\ldots,n_{\lfloor r/s\rfloor},
\underbrace{0,\ldots,0}_{s-1+r-\lfloor r/s\rfloor s})\\
&+&B(n_1,b_1n_1,\ldots,b_{s-2}n_1,n_1,n_2,b_1n_2,\ldots,b_{s-2}n_2,n_2,\ldots,\\
&& \quad n_{\lfloor r/s\rfloor},
b_1n_{\lfloor r/s\rfloor},\ldots,b_{s-2}n_{\lfloor r/s\rfloor},n_{\lfloor r/s\rfloor},
\underbrace{0,\ldots,0}_{r-\lfloor r/s\rfloor s})\\
&=&1
\end{eqnarray*}
for $(n_1,n_2,\ldots,n_{\lfloor r/s\rfloor})\not=(0,0,\ldots,0)$.
As before, the polynomial
\begin{eqnarray*}&&F(X_1,X_2,\ldots,X_{\lfloor r/s\rfloor})\\
&=&B(X_1,0,\ldots,0,X_2,0,\ldots,0,...)\\&&+B(X_1,b_1X_1,\ldots,b_{s-2}X_1,X_1,X_2,b_1X_2,\ldots,b_{s-2}X_2,X_2,\ldots)
\end{eqnarray*}
has degree $\lfloor r/s\rfloor$ and sparsity $2^{\lfloor r/s\rfloor}-1$.
The result follows from $\deg(B)\ge \deg(F)$ and ${\rm spr}(B)\ge \lceil \frac{{\rm spr}(F)}{2}\rceil$.\hfill $\Box$\\

Taking a random (sufficiently large) prime, with very high probability $>0.999$ we can take $s\le 5$, that is ${\rm N}(p)\le 31$, and get
$$\deg(B)\ge 0.2 \log_2 p +O(1)$$
which improves $(\ref{degigor})$.
With probability $>0.9$ we still have $s\le 3$, that is ${\rm N}(p)\le 7$, and 
$${\rm spr}(B)\gg p^{1/3}$$
which improves $(\ref{sprigor})$.

Numerical data suggests that both the degree and the sparsity of $B$ are close to the expected values, in particular, 
$$\deg(B) \ge r-2\quad\mbox{for}\quad 2 < p < 10000.$$
Moreover, Figure \ref{fig:sparsity} in Appendix 1 supports the conjecture that 
$${\rm spr}(B) = 2^{r-1} + O(p^{1/2}).$$



\section{Bounds on the $N$th linear complexity}\label{seq}

In this section we prove bounds on the $N$th linear complexity of sequences with the property $s_{2n}=1-s_n$, $n=1,2,\ldots$ which includes the sequences corresponding to the Liouville functions. In particular, we substantially improve results on the Legendre sequence of period $p$ for $p\equiv \pm 3\bmod 8$. For $p\equiv \pm 1\bmod 8$ the lower bound
$(\ref{leglow})$ is currently the best known one.

The proof of the general bound uses a result on the $N$th lattice level defined below which is closely related to the $N$th linear complexity, see the following subsection.

For the Legendre sequence of period $p\equiv \pm 1\bmod 8$ we provide some numerical data which leads to the conjecture 
$$L({\cal L}_p,N)=\frac{N}{2}+O(\log p),\quad N=1,2,\ldots,p+1.$$

\subsection{$N$th lattice level}

We recall a measure of pseudorandomness closely related to the $N$th linear complexity.
A binary sequence ${\cal S}=(s_n)_{n=1}^\infty$ passes the
$S$-dimensional $N$-lattice test if the vectors
$$\{(s_n-s_1,s_{n+1}-s_2,\ldots,s_{n+S-1}-s_S): n=2,3,\ldots,N-S+1\}$$
span $\mathbb{F}_2^S$.
The greatest $S=S({\cal S},N)$ such that ${\cal S}$ passes the $S$-dimensional $N$-lattice test is called the {\em $N$th lattice level} of ${\cal S}$.
By \cite[Proposition~4]{DW03} we have $S({\cal S},N)\le \lfloor N/2\rfloor$
and those sequences which attain this bound for all $N$ are characterized in the following proposition.

\begin{proposition}
\label{prop1}
\emph{\cite[Theorem~22]{DMW04}}\\
The sequence ${\cal S}=(s_n)_{n=1}^\infty$ satisfies
$$S({\cal S},N)=\left\lfloor \frac{N}{2}\right\rfloor$$
if and only if
$$s_{2n}=1-s_n,\quad n=1,2,\ldots$$
\end{proposition}

We have the following strong connection between the $N$th linear complexity and the $N$th lattice level.

\begin{proposition}
\label{prop2}
\emph{\cite[Theorem~1]{DW03}}\\
We have either
$$S({\cal S},N)=\min\{L({\cal S},N),N+1-L({\cal S},N)\}$$
or
$$S({\cal S},N)=\min\{L({\cal S},N),N+1-L({\cal S},N)\}-1.$$
\end{proposition}

\subsection{Linear complexity}

By the following result we can obtain upper bounds on the $N$th linear complexity from suitable lower bounds. 
\begin{proposition}
\label{prop3}
\emph{\cite[Lemma~5]{DW04}}\\
Let ${\cal U}=(u_n)_{n=1}^\infty$ be a sequence with $u_1\le 0$ and $u_n\le u_{n-1}+1$ for $n\ge 2$.\\
If the sequence ${\cal S}$ satisfies
$$L({\cal S},N)\ge u_N \quad \mbox{for }N\ge 2,$$
then we have
$$L({\cal S},N)\le N-u_{N-1}.$$
\end{proposition}

Now we are able to prove results on the $N$th linear complexity of many sequences including those corresponding to the Liouville functions.

\begin{corollary}
\label{cor1}
  If the sequence ${\cal S}=(s_n)_{n=1}^\infty$ satisfies
  $$s_{2n}=1-s_n,\quad n=1,2,\ldots$$
  then 
  $$
  \left\lfloor \frac{N}{2}\right\rfloor
  \le L({\cal S},N)\le \left\lfloor \frac{N}{2}\right\rfloor+1,\quad  N=1,2,\ldots$$ 
\end{corollary}
Proof. 
Combining Propositions \ref{prop1} and \ref{prop2} we get the lower bound.
Taking $$u_n=\left\lfloor\frac{n}{2}\right\rfloor\quad \mbox{for }n\ge 1$$ 
we get the upper bound by Proposition~\ref{prop3}. \hfill $\Box$\\

We can also adjust this result to the Legendre sequence in the case when $2$ is a quadratic non-residue modulo $p$.
\begin{theorem}
Let $p\equiv \pm 3 \bmod 8$ be a prime and ${\cal L}_p$ the $p$-periodic Legendre sequence.
Then we have 
$$ \left\lceil \frac{\min\{N,2p-1\}-1}{2}\right\rceil \le L({\cal L}_p,N)\le
\left\lfloor\frac{\min\{N,2p-2\}}{2}\right\rfloor+1.$$
\end{theorem}
Proof. Consider the sequence ${\cal S}=(s_n)_{n=1}^\infty$ with
$$s_n=\ell_n\quad \mbox{if }  n\not\equiv 0\bmod 2p$$ 
and 
$$s_{2kp}=1-s_{kp}\quad \mbox{for }k=1,2,\ldots$$
${\cal S}$ satisfies the conditions of Corollary~\ref{cor1}.
Now the first $2p-1$ elements of ${\cal S}$ and~${\cal L}_p$ coincide and we have
$$L({\cal L}_p,N)=L({\cal S},N),\quad N=1,2,\ldots,2p-1,$$
and the result follows. 
For $N\ge 2p$ we have $L({\cal L}_p,N)=L({\cal L}_p)\in \{p-1,p\}$
by~$(\ref{leglincompl})$. \hfill $\Box$\\

For the case $p\equiv \pm 1\bmod 8$,
Figure \ref{fig:lc1} and Figure \ref{fig:lc2} in Appendix 1
support the conjecture that 
$$L({\cal L}_p,N)=\frac{N}{2}+O(\log p),\quad N=1,2,\ldots,p+1$$
and 
$$\max_{N=1,2,\ldots,p+1}\left|L({\cal L}_p,N)-\frac{N}{2}\right|$$
is of order of magnitude $\log p$.

\section{Conclusion}
We improved results on the degree and sparsity of Boolean functions representing the Legendre symbol
and the $N$th linear complexity of the Legendre sequence. We presented these result in a more general form for all arithmetic functions $f$ with the property $f(2n)=-f(n)$, $n=1,2,\ldots$. For example, besides the Legendre symbol this includes also the Boolean function representing both the integer and the polynomial Liouville function. 

\section*{Acknowledgment}
The authors wish to thank Zhixiong Chen and the reviewers for several useful comments.

\section*{Appendix 1: Figures}

\begin{figure}[H]
\caption{ The distance of ${\rm spr}(B)$ from $2^{r-1}$ for all primes $2 < p < 10000$, where $B$ is a Boolean function corresponding to the Legendre symbol with modulus $p$ and $r = \lfloor\log_2 p\rfloor$.}
\label{fig:sparsity}
\centering
\includegraphics[width=\textwidth]{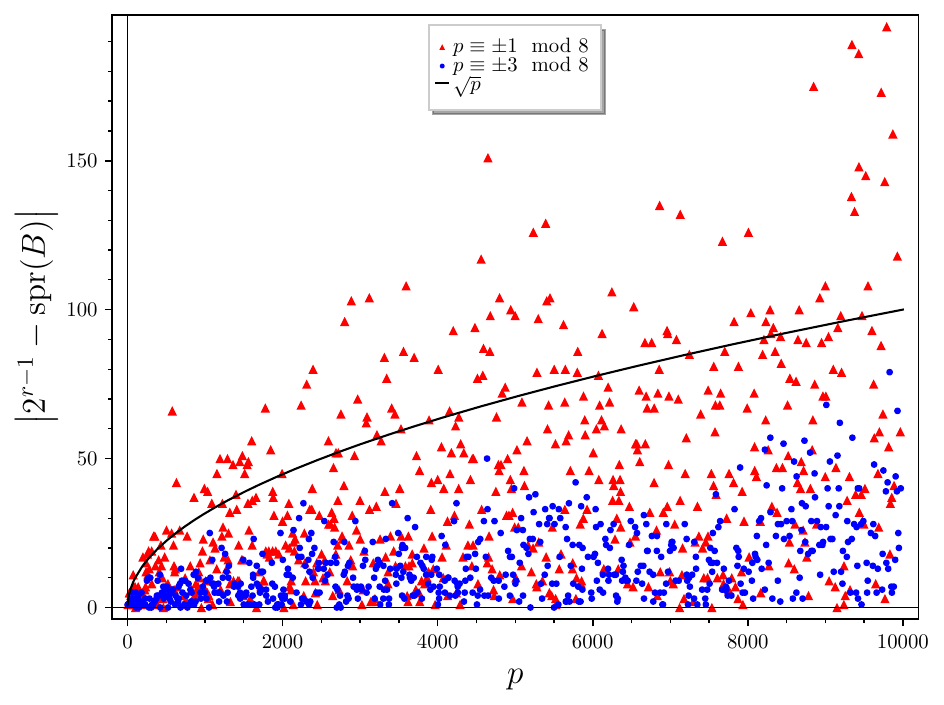}
\end{figure}

\begin{figure}[H]
\caption{ The maximum distance of $L(\mathcal{L}_p, N)$ from $N/2$, $N = 1, \ldots, p + 1$, for all primes $p < 10000$ with $p\equiv \pm 1\bmod 8$.}
\label{fig:lc1}
\centering
\includegraphics[width=\textwidth]{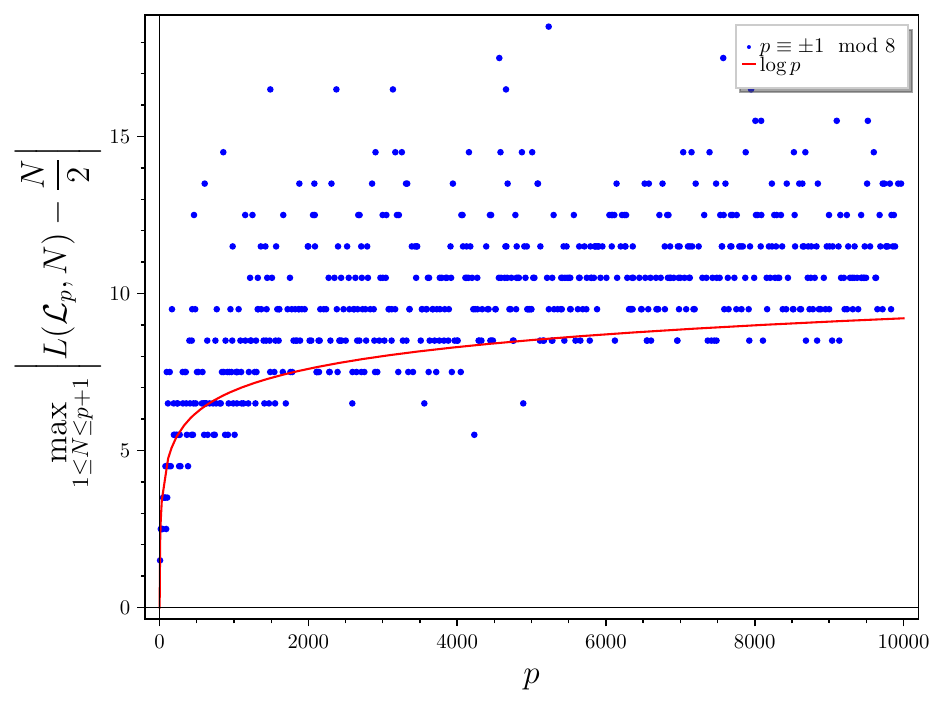}
\end{figure}

\begin{figure}[H]
\caption{$L(\mathcal{L}_{100049}, N) - \frac{N}{2}$ for $N = 1, \ldots, 100050$.}
\label{fig:lc2}
\centering
\includegraphics[width=\textwidth]{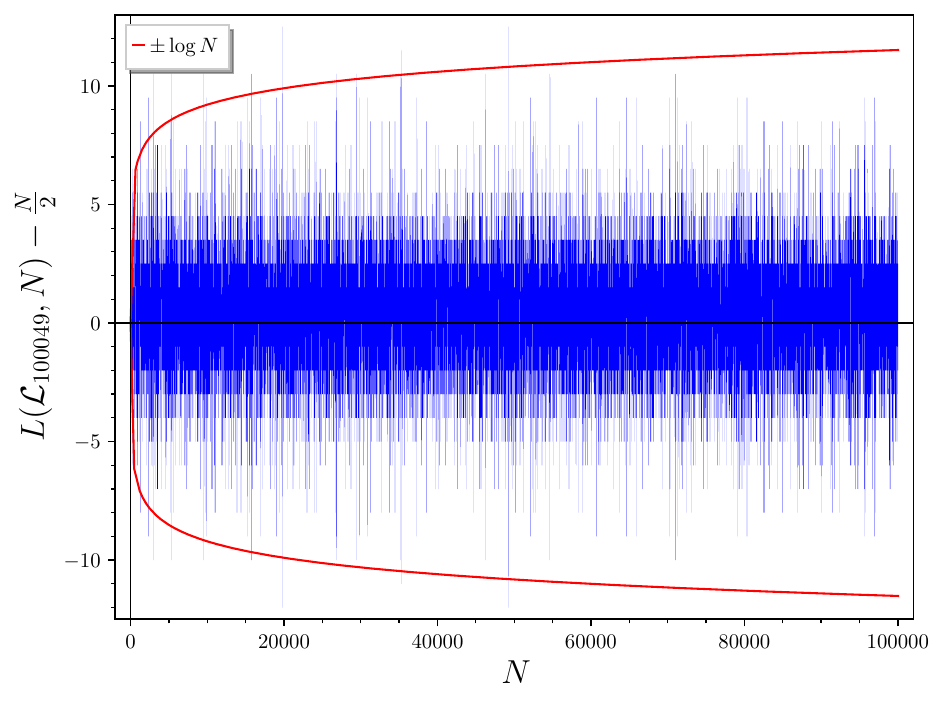}
\end{figure}

\section*{Appendix 2:\\ Correlation measure and linear complexity}
The {\em $N$th correlation measure $C_k({\cal S},N)$ of order $k$} of a binary sequence ${\cal S}=(s_n)_{n=1}^\infty$ was introduced by Mauduit and S\'ark\"ozy in \cite{MS97},
$$C_k({\cal S},N) = \max_{M,D}\left|\sum_{n=1}^M(-1)^{s_{n+d_1}+s_{n+d_2}+\ldots+s_{n+d_k}}\right|,$$
where the maximum is taken over all integer vectors $D = (d_1,\ldots,d_k)$ and integers $M$ such that
$0\le d_1<d_2<\ldots<d_k\le N-M$. \\

Improving a relation of Brandst\"atter and the second author \cite{brwi06} between linear complexity and correlation measure, Gomez et al. \cite[Corollary~4]{CGGT22}
proved the following result:\\
Let $K$ and $N$ be positive integers with $2 \le K^2 < N$. If a binary sequence ${\cal S}$
satisfies $C_k({\cal S}, N ) < N/2$ for every $k \le K$, then we have
$$L(S, N ) \gg K \log(N),$$
where the implied constant is absolute.

For example, for the Legendre sequence ${\cal L}_p$ we have for $N\le p$,
$$C_k({\cal L}_p,N)=O(kp^{1/2}\log p),$$
see $\cite{MS97}$, and thus
$$L({\cal L}_p,N)\gg \frac{N\log(N)}{p^{1/2}\log(p)}.$$
Since otherwise the result is trivial we may assume $N\gg p^{1/2}$ and thus
$$L({\cal L}_p,N)\gg \frac{N}{p^{1/2}},$$
which improves the bound of $\cite{brwi06}$ by a factor $\log(p)$.

\end{document}